\documentclass[12pt]{amsart}
\usepackage[utf8]{inputenc}
\usepackage{ulem}
\usepackage{subfig}
\usepackage{lipsum, setspace}
\allowdisplaybreaks[1] 
\usepackage{amssymb,amsmath,amsthm, graphicx, epstopdf}

\usepackage{enumerate}
\usepackage{xcolor}
\usepackage{blkarray}
\usepackage{mathrsfs}
\newtheorem{theorem}{Theorem}[section]
\newtheorem{proposition}[theorem]{Proposition}
\newtheorem{lemma}[theorem]{Lemma}

\theoremstyle{definition}
\newtheorem{definition}[theorem]{Definition}
\newtheorem{example}[theorem]{Example}
\newtheorem{remark}[theorem]{Remark}

\newcommand{\clh}{\mathcal{H}}

\DeclareRobustCommand{\rrho}{{\mathpalette\irrho\relax}}
\newcommand{\irrho}[2]{\raisebox{\depth}{$#1\rho$}}

\usepackage[utf8]{inputenc}
\newcommand\hlight[1]{\tikz[overlay, remember picture,baseline=-\the\dimexpr\fontdimen22\textfont2\relax]\node[rectangle,fill=white!50,rounded corners,fill opacity = 0.2,draw,thick,text opacity =1] {$#1$};}

\usepackage[bookmarksnumbered, colorlinks, plainpages]{hyperref}
\title[Pseudo $S$-spectra of quaternionic operators]{Pseudo $S$-spectra of special operators in quaternionic Hilbert spaces} 

\author[Dhara]{Kousik Dhara}
\address{Kousik Dhara, Department of Mathematics, 
 Weizmann Institute of Science, Rehovot 7610001,
Israel}
\email{kousik.dhara@weizmann.ac.il}

\author[Pamula]{Santhosh Kumar Pamula}
\address{Santhosh Kumar Pamula,
Department of Mathematical Sciences,
Indian Institute of Science Education and Research (IISER) Mohali, 
SAS Nagar, Manauli Post 140306, Punjab, India}
\email{santhoshkp@iisermohali.ac.in}


\numberwithin{equation}{section}

\begin{document}
\maketitle
\begin{abstract}
For a bounded quaternionic operator $T$ on a right quaternionic Hilbert space $\mathcal{H}$  and $\varepsilon >0$, the pseudo $S$-spectrum of $T$ is defined as
\begin{align*}
\Lambda_{\varepsilon}^{S}(T) := \sigma_S (T) \bigcup \left \{ q \in \mathbb{H}\setminus \sigma_S(T):\; \|Q_q(T)^{-1}\| \geq \frac{1}{\varepsilon} \right\},
\end{align*}
where $\mathbb{H}$ denotes the division ring of quaternions, $\sigma_S(T)$ is the $S$-spectrum of $T$ and $Q_{q}(T)= T^2-2 \text{Re}(q)T+|q|^2I$. This is a natural generalization of pseudospectrum from the theory of  complex Hilbert spaces.
In this article, we investigate several properties of the pseudo $S$-spectrum and explicitly compute the pseudo $S$-spectra for  some special classes of operators such as upper triangular matrices, self adjoint-operators, normal operators and orthogonal projections.   In particular, by an application of $S$-functional calculus, we show that a quaternionic operator is a left multiplication operator induced by a real number $r$  if and only if for every $\varepsilon>0$ the pseudo $S$-spectrum of the operator is the circularization of a closed disc in the complex plane centered at $r$ with the radius $\sqrt{\varepsilon}$.  Further, we propose a $G_1$-condition for quaternionic operators and prove some results in this setting. 
\end{abstract}
\section{Introduction}
In the theory of complex Hilbert spaces, pseudospectrum plays a vital role in understanding several properties of a non self-adjoint operator.  Formally, for $\varepsilon>0$, the $\varepsilon$-pseudospectrum (in short, pseudospectrum) of a bounded linear operator $A$ on a complex Hilbert space $\mathcal{K}$ is defined by
\begin{align*}
    \Lambda_{\varepsilon}(A) = \sigma (A) \cup \left \{ \lambda \in \mathbb{C}\setminus \sigma(A):\; \|(A-\lambda I)^{-1}\| \geq \frac{1}{\varepsilon} \right\}.
\end{align*}
Here the spectrum $\sigma(A)$ of $A$  consists of all complex numbers $\lambda$ for which $(A-\lambda I)$ does not have a bounded inverse on $\mathcal{K}$.  In particular, if $A$ is a normal operator, then the  $\varepsilon$-pseudospectrum of $A$ coincides with the closed $\varepsilon$-disc around the spectrum of $A$.

Originally, the notion of pseudospectra arose from several observations such as the existence of $\varepsilon$-approximate eigenvalues that are far from  the spectrum \cite{Landau2}, the instability of spectrum under slight perturbations and the necessity to estimate the resolvent norm of an operator.
The substantial monograph by L. N. Trefthen and M. Embree \cite{Tref2} provides a detailed overview of the development of pseudospectral theory and its applications to diverse fields. For instance, it is extensively used in the study of spectral approximation of Toeplitz matrices, random matrices, stability of numerical solutions of differential equations and control theory. 
We refer to H. J. Landau \cite{Landau2,Landau1}, L. N. Trefthen \cite{Tref1}, A. Böttcher \cite{Bottcher1,Bottcher2}, E. B. Davies \cite{Davies}  for some classical work on pseudospectra. Also see \cite{Cui1,Cui2,Dhara2,Dhara, Frommer, Krishnan,Seidel,Shargo} for recent developments on pseudospectral theory.

In 1936, Birkhoff and von Neumann showed that quantum mechanics can be formulated  on the division ring of quaternions along with real and complex fields \cite{Birkhoff}. Since then, several researchers have studied quaternionic Hilbert spaces and made significant development in this direction (for example, see \cite{Fink, Viswanath} and references therein). They used the notion of left and right spectrum, but these spectra were insufficient to construct the quaternionic spectral theory. Only in 2006, F. Colombo and I. Sabadini  introduced the notion of $S$-spectrum which led to the development of spectral theory in quaternionic Hilbert spaces. We refer the reader to Subsection 1.2.1 of the book \cite{Colombo2} by F. Colombo, J. Gantner and D.P. Kimsey for a brief discussion and motivation of the $S$-spectrum. In this book \cite{Colombo2}, the authors described $S$-functional calculus in the introduction and also explained how Hypercomplex analysis methods were used to identify the appropriate notion of the quaternionic spectrum whose existence was suggested by quaternionic quantum mechanics.  In 2016, D. Alpay, F. Colombo and D.P. Kimsey gave a full proof of the spectral theorem for both bounded and unbounded quaternionic operators (see \cite{Alpay-0} for details).

In fact, the quaternionic spectral theory based on the $S$-spectrum, that is spread over several papers, has been systematically organized, with historical notes, in the books \cite{Colombo3,Colombo2}. For the spectral theory on  S-spectrum for Clifford operators, we refer the reader to the book \cite{Colombo1} and the references
therein. The seminal work on the $S$-spectrum and the full development of the S-functional calculus are contained in the papers \cite{Alpay-R,Alpay1, Alpay2, Colombo-S2,Colombo-S1}.  Additionally, we highlight the detailed description of the main research directions in the spectral theory based on the $S$-spectrum as follows: for perturbation of quaternionic normal operators, see \cite{Paula};  for the intrinsic S-Functional Calculus and Spectral Operators, see \cite{Gantner}; for Slice Hyperholomorphic Schur Analysis, see \cite{Alpay-S2}; for Quaternionic de Branges spaces and characteristic operator function, see \cite{Alpay-S1} and for a recent overview and applications of spectral theory, see \cite{Colombo-P}. Also, for a multiplication form of the quaternionic normal operator based on the $S$-spectrum, see \cite{Golla-Pamula 2020};  for Borel functional calculus for quaternionic normal operators, see \cite{Ramesh} and for a brief discussion of representation of quaternionic compact operators, see \cite{Golla-Pamula 2017}.

For a bounded right linear operator $T$ on a quaternionic Hilbert space $\mathcal{H}$ ($T \in \mathcal{B}(\mathcal{H})$ for short), the \textit{$S$-spectrum} of $T$ is defined by
\begin{align*}
\sigma_S(T)=\left \{ q\in \mathbb{H}:Q_{q}(T) \text{ is not invertible in } \mathcal{B}(\mathcal{H}) \right \},
\end{align*}
where $Q_{q}(T)= T^2- 2 \text{Re}(q)T+|q|^2I$. Here Re($q$) and $|q|$ denote the the real part of the quaternion $q$ and the modulus of $q$ respectively. The notion of $S$-spectrum is a natural extension of the concept of right eigenvalues for quaternion matrices to infinite dimensional Hilbert spaces. 

As the theory of pseudospectra progresses, it is natural to study this notion in the quaternionic setting with possible applications.
As in the classical case, the motivation comes from the operator equation of the form  $Q_{q}(T)x=y$ with the assumption that $q\notin \sigma_S(T)$. The study of solutions that remain stable under perturbation of $y$ is of general interest. Suppose that $Q_{q}(T)x' = y'$ is the perturbed equation  with $\|y-y' \|< \varepsilon$. Then we have \begin{equation*}
    \|x-x'\|<\varepsilon \|Q_{q}(T)^{-1}\|.
\end{equation*}
This leads to the requirement of  estimating norm of the pseudo-resolvent operator $\|Q_{q}(T)^{-1} \|$ and we need to ensure that $\varepsilon$ is small compared to $\|Q_{q}(T)^{-1} \|$. It gives rise to the following definition:
For $T \in \mathcal{B}(\mathcal{H})$ and $\varepsilon >0$, the \textit{$\varepsilon$-pseudo $S$-spectrum} (in short, pseudo $S$-spectrum) of $T$ is defined by 
\begin{align} \label{Definition: pseudo}
\Lambda_{\varepsilon}^{S}(T) = \sigma_S (T) \cup \Big\{ q \in \mathbb{H}\setminus \sigma_S(T):\; \big\|Q_q(T)^{-1}\big\| \geq \frac{1}{\varepsilon} \Big\}.
\end{align}

Here, it is worth to note that in contrast to the complex case, one faces two challenges in the quaternionic setup: the non-commutative nature of  the quaternions and non-invertibility of a second order bounded linear operator. This makes the study of pseudo $S$-spectrum to be intriguing in the present framework. 

The notion of pseudo $S$-spectrum possibly first appeared in the recent work of  A. Ammar, A. Jeribi and N. Lazrag in \cite{Ammar}. The authors studied some basic properties and gave a characterization for the Weyl pseudo $S$-sepctrum in terms of Fredholm operators.

In this paper, we systematically study the pseudo $S$-spectrum of a quaternionic operator, prove elementary properties and often compare the obtained results with the complex set up. The computation of pseudo $S$-spectra for quaternionic operators is not found in the literature and we believe that it is necessary to provide some examples in this aspect {(for instance, see Section 4 of this article).} In this direction,  we obtain a concrete description of pseudo $S$-spectra of some special classes of operators including $(2\times 2)$-upper triangular matrices, self adjoint operators, normal operators and orthogonal projections.
For instance, we establish that the pseudo $S$-spectrum of a self adjoint operator is the union of the circularizations of the discs centred at $r$ and radius $\sqrt{\varepsilon}$ with $r$ runs over the $S$-spectrum (see Theorem \ref{Self adjoint S-spectrum}). 
In particular, a quaternionic operator is a left multiplication operator induced by a real number $r$  if and only if for every $\varepsilon>0$, the pseudo $S$-spectrum of the operator is the circularization of a closed disc in the complex plane which is centered at $r$ and radius $\sqrt{\varepsilon}$ (see Proposition \ref{Prop: Mult op real}). In fact, a detailed investigation is done on the pseudo $S$-spectrum of a left multiplication operator induced by a quaternion (see Theorems \ref{bound for s-pseudo} and \ref{Theorem: emptyorconnected}). We wish to emphasize that some of the techniques used in this paper rely on the the notion of circularization and $S$-functional calculus (see \cite{Alpay-R, Colombo-S2, Colombo-S1, Colombo1}) for quaternionic operators.

We organize the article as follows: In Section 2, we fix some notations, review some basic definitions and known results that are useful for the following sections. Section 3 is dedicated to describe some properties of pseudo $S$-spectrum. In Section 4, we compute the pseudo $S$-spectra of  special operators and propose a first order growth condition (or $G_{1}$-condition) for quaternionic operators and study their properties.

\section{Notation and Preliminaries}
The set of all quaternions with real coefficients is denoted by $\mathbb{H}$ and it is formally defined as
\begin{equation*}
\mathbb{H} = \big\{  q_{0}+ q_{1}i+q_{2}j+q_{3}k:\; q_{\ell} \in \mathbb{R},\; \ell = 0,1,2,3\big\}.
\end{equation*}
Here $i,j$ and $k$ are known as the unit imaginary quaternions satisfying
\begin{equation}\label{Equation: Multiplication}
i^{2} = j^{2} = k^{2} = ijk= -1.
\end{equation}
The above relation directs the multiplication in $\mathbb{H}$. Especially, we see that 
\begin{equation*}
ij = -ji=k;\; ki=-ik=j;\; jk=-kj=i.
\end{equation*}
We identify $a\in \mathbb{R}$ with $a 1$ in $\mathbb{H}$ and $\alpha + \beta i \in \mathbb{C}$ with $\alpha + \beta i + 0 j + 0 k$ and hence $\mathbb{H}$ is a noncommutative real division algebra (contains $\mathbb{R}$ and $\mathbb{C}$)  with respect to real vector space structure and the mutliplication given by the Equation (\ref{Equation: Multiplication}). For $q \in \mathbb{H}$, the real and imaginary part of $q$ are given by $\text{Re}(q) = q_{0}$ and $\text{Im}(q) = q_{1}i+q_{2}j+q_{3}k $ respectively. In fact, for a subset $S$ of $\mathbb{H}$, we write $\text{Re}(S) = \{ \text{Re}(s):\; s \in S  \}$ and $\text{Im}(S) = \{\text{Im}(s):\; s\in S\}$. The modulus of  $q$ is given by $|q| = \sqrt{q_{0}^{2}+q_{1}^{2}+q_{2}^{2}+q_{3}^{2}}$. The unit sphere of purely imaginary quaternions is defined by
\begin{equation*}
    \mathbb{S} = \Big\{q\in \mathbb{H}:\; \overline{q} = -q\; \text{and}\; |q|=1\Big\}.
\end{equation*}
Note that $m\in \mathbb{S}$ if and only if $m^{2}=-1$, so the elements of $\mathbb{S}$ are often regarded as imaginary units. In fact, $\mathbb{S}$ is a two dimensional sphere in $\mathbb{R}^{4}.$ For each $m \in \mathbb{S}$, the real algebra generated by $1$ and $m$ is denoted by $\mathbb{C}_{m} = \{\alpha + \beta m:\; \alpha, \beta \in \mathbb{R}\}.$ Notice that each $\mathbb{C}_{m}$ is isomorphic to the complex plane and $\mathbb{C}_{m} \cap \mathbb{C}_{n} = \mathbb{R}$ for $m\neq \pm n.$ For every  $q \in \mathbb{H}\setminus \mathbb{R}$, we associate a imaginary unit $m_{q}:= \frac{\text{Im}(q)}{|\text{Im}(q)|} \in \mathbb{S}$ in the sense that $q= \text{Re}(q)+ |\text{Im}(q)|\; m_{q} \in \mathbb{C}_{m_{q}}$ and hence $\mathbb{H} = \bigcup\limits_{m \in \mathbb{S}}\mathbb{C}_{m}.$ Now  $p,q \in \mathbb{H}$ are said to be {\it similar}  if $p = s^{-1}qs$ for some  $s \in \mathbb{H}\setminus \{0\}.$ 
This is an equivalence relation on $\mathbb{H}$ and the equivalence class of $q$ is given by 
\begin{equation*}
    [q] = \big\{ \text{Re}(q)+ m |\text{Im}(q)|:\; m \in \mathbb{S}\big\}.
\end{equation*}
In view of this, we have the following definition.
\begin{definition} Let $X$ be a nonempty subset of $\mathbb{C}$.  The \textit{circularization} of $X$ is defined by
\begin{equation*}
    \Omega_{X} = \big\{ \alpha + \beta m:\; \alpha + \beta i \in X, m \in \mathbb{S},\alpha, \beta \in \mathbb{R}\big\} \subseteq \mathbb{H}.
\end{equation*}
In other words, $\Omega_{X} = \bigcup\limits_{x\in X} [x].$ Moreover, a nonempty subset ${K}$ of $\mathbb{H}$ is said to be {\it circular} or {\it axially symmetric} if $K = \Omega_{X}$ for some nonempty subset $X$ of $\mathbb{C}.$ 
 \end{definition}
 For example, the circularization of $i$ is $\Omega_{\{i\}} = \mathbb{S}$. For $\lambda \in \mathbb{C}$ and $r>0$, we denote    the open disk with center $z$ and radius $r$ by $D(\lambda, r) = \big\{z\in \mathbb{C}:\; |z-\lambda|<r\big\}$ . Then $\Omega_{D(\lambda,\; r)}$ is a solid torus around $X$-axis and whose intersection with the complex plane is $D(\lambda,\; r) \cup D(\overline{\lambda}, r)$. Now we recall the definition of a quaternionic Hilbert space. 

\begin{definition}
Let $\mathcal{H}$ be a right $\mathbb{H}$-module. A map $\langle \cdot, \cdot \rangle\colon \mathcal{H} \times \mathcal{H} \to \mathbb{H}$ is said to be an inner product on $\mathcal{H}$ if it satisfies the following conditions:
\begin{enumerate}
\item {\bf Positivity}: $\langle x, x\rangle \geq 0$. Moreover, $\langle x, x\rangle = 0$ if and only if $x=0$
\item {\bf Right linearity}: $\langle x, y q+ z\rangle = \langle x, y\rangle q + \langle x, z\rangle$
\item {\bf Quaternionic Hermiticity}:$\langle x, y \rangle = \overline{\langle y, x\rangle},$
\end{enumerate}
for every $x,y,z\in \mathcal{H}$ and $ q \in \mathbb{H}$. Furthermore, if $\mathcal{H}$ is complete with respect to the induced norm given by $\| x\|: = \langle x, x\rangle^{\frac{1}{2}}$, then $\mathcal{H}$ is said to be a {\it{right quaternionic Hilbert space}}. 
\end{definition}

 Throughtout this article, $\mathcal{H}$ denotes a right quaternionic Hilbert space. A map $T \colon \mathcal{H} \to \mathcal{H}$ is said to be  {\it right $\mathbb{H}$-linear} or {\it quaternionic operator} if 
\begin{align*}
T(xq+y)=T(x)q+T(y) \, \, \, \text{ for all } x,y\in \mathcal{H} \text{ and } q\in \mathbb{H}.
\end{align*}
Further, $T$ is said to be {\it bounded}  if there is $M>0$ such that $\|Tx\| \leq M \|x\|$ for all $x \in \mathcal{H}$.  
We denote the space of all bounded quaternionic operators on $\mathcal{H}$ by $\mathcal{B}(\mathcal{H})$. Throughout, by a quaternionic operator, we mean a bounded quaternionic operator. Note that $\mathcal{B}(\mathcal{H})$ is  a real Banach algebra equipped with the operator norm given by 
\begin{align*}
\|T\|= \sup\big\{\|Tx\|:\; x\in \mathcal{H}, \|x\|\leq 1\big\}.
\end{align*}  For every $T \in \mathcal{B}(\mathcal{H})$ there is a unique operator $T^{\ast} \in \mathcal{B}(\mathcal{H})$, called the {\it adjoint} of $T$, satisfying
\begin{equation*}
    \langle x, Ty\rangle = \langle T^{\ast}x,\; y\rangle\; \text{for all} \; x,y \in \mathcal{H}.
\end{equation*}
An operator $T \in \mathcal{B}(\mathcal{H})$ is said to be: (i) {\it self-adjoint} if $T^{\ast} = T$; (ii) {\it normal} if $TT^{\ast} = T^{\ast}T$; (iii) {\it orthognal projection} if $T^{\ast}=T=T^2$; (iv) {\it hyponormal} if $T^{\ast}T-TT^{\ast} \geq 0$; and (v) {\it subnormal} if there exist a normal operator $S$ on some quaternionic Hilbert space $\widehat{\mathcal{H}}$ such that $\mathcal{H}$ is a closed subspace of $\widehat{\mathcal{H}}$, $S(\mathcal{H})\subseteq \mathcal{H}$ and $S\vert_{\mathcal{H}} = T.$

The $S$-spectrum $\sigma_{S}(T)$ of $T$ is a non-empty compact subset of $\mathbb{H}$. The \textit{$S$-resolvent} of $T$ is defined as  $\rrho_{S}(T) = \mathbb{H}\setminus \sigma_{S}(T)$ and the {\it $S$-spectral radius} is defined by 
\begin{align*}
r_S(T)=\max \{|q|: q\in \sigma_S(T) \}.
\end{align*}
Also the $S$-spectral radius formula is given by $r_S(T)=\underset{n\rightarrow \infty}{\lim} \|T^n\|^{1/n}$.
We refer to the classic book \cite{Colombo2} by Colombo, Ganter and Kimsey for a comprehensive overview of the study of $S$-spectrum and its historical roots in quaternionic Hilbert spaces.

\section{Properties of pseudo $S$-spectrum }
In this section, we study some elementary properties of pseudo $S$-spectrum of a quaternionic operator. 
 The following proposition gives some equivalent definitions of pseudo $S$-spectrum. 
\begin{proposition}\label{Proposition: singular}
Let $T\in \mathcal{B}(\mathcal{H})$ and $\varepsilon>0$. Then the following statements are equivalent:
\begin{enumerate}
\item $q\in \Lambda_\varepsilon^S(T).$
\item Either $q\in \sigma_S(T)$ or there exists $u\in \mathcal{H}$ with $\|u\|=1$ such that $\|Q_{q}(T)u\|\leq \varepsilon. $
\end{enumerate}
Moreover, if $\mathcal{H}$ is finite dimensional, then
\begin{equation*}
    \Lambda_\varepsilon^S(T)=
\{q\in \mathbb{H}: s_\text{min}(Q_{q}(T))\leq \varepsilon \},
\end{equation*}
where $s_\text{min}(Q_{q}(T))$ is the smallest singular value of $Q_{q}(T)$.
\end{proposition}
\begin{proof} The proof follows similar lines as in the complex case (for instance, see \cite{Tref2}). 
\end{proof}

In the following proposition, we collect some elementary properties of pseudo $S$-spectrum. Some of these can be compared to the existing classical results from \cite{Tref2}.  
\begin{proposition}\label{Elementary prop}
Let $T \in \mathcal{B}(\mathcal{H})$ and $\varepsilon>0$. Then the following statements hold:
\begin{enumerate}
\item $\Lambda_\varepsilon^S(T)$ is axially symmetric. 
\item $ \Lambda_{\varepsilon_1}^{S}(T)\subseteq \Lambda_{\varepsilon_2}^{S}(T) \; \text{whenever}\; 
0< \varepsilon_1 < \varepsilon_2$.
\item $\sigma_S (T)=\underset {\varepsilon>0}\bigcap \Lambda_{\varepsilon}^{S}(T).$
\item $\Lambda_{\varepsilon}^{S}(r T)= r \cdot \Lambda_{\frac \varepsilon {r^2}}^{S}(T)$ for all $ r \in \mathbb{R}\setminus \{0\}$.

\item Suppose  $V\in \mathcal{B}(\mathcal{H})$ is invertible and if $B=V^{-1}TV$, then
\begin{align*}
\Lambda_{\frac{\varepsilon}{k}}^S(T) \subseteq \Lambda_\varepsilon^S(B)\subseteq \Lambda_{\varepsilon k}^S(T),
\end{align*}
where $k=\|V^{-1}\|\|V\|$.

\item  $\Lambda_\varepsilon^S(T_{1}\oplus T_{2})= \Lambda_\varepsilon^S(T_{1})
\bigcup \Lambda_\varepsilon^S(T_{2})$ for $T_{1}, T_{2} \in \mathcal{B}(\mathcal{H})$.

\item $ \Lambda_{\varepsilon}^{S}(T)\subseteq \{ q\in \mathbb{H}: |q|\leq \|T\|+\sqrt{\varepsilon}\}$. 
In particular, $\Lambda_{\varepsilon}^{S}(T)$ is a non-empty compact subset of $\mathbb{H}$.
\end{enumerate}
\end{proposition}

\begin{proof} 
Since $Q_{p}(T) = Q_q(T)$ whenever $p \in [q]$, it implies that $\Lambda_{\varepsilon}^{S}(T)$ is axially symmetric. The assertions (2) and (3) follow directly from the Definition \ref{Definition: pseudo}.

\noindent Proof of (4). For $q \in \mathbb{H}$ and $r \in \mathbb{R}\setminus \{0\}$, we see that 
\begin{equation}\label{Equation: rq}
    Q_q(rT)=r^2 \left(T^2 - 2 \text{Re}\left(\frac{q}{r}\right)T + \left|\frac{q}{r} \right|^2 \right)
= r ^2 Q_{\frac{q}{r}}(T).
\end{equation}
Then  $\sigma_S(rT)=r\cdot \sigma_S(T)$. Now for $q \in \Lambda^{S}_{\varepsilon}(rT)\setminus \sigma_{S}(rT)$, it follows that 
\begin{align*}
    \|Q_{q}(rT)^{-1}\| \geq \frac{1}{\varepsilon} \; \iff \; \|Q_q(T)^{-1}\| \geq \frac{1}{(\frac{\varepsilon}{r^{2}})}.
\end{align*}

\noindent Proof of (5). First note that $\sigma_{S}(B) = \sigma_{S}(T)$, since $Q_{p}(B) = V^{-1}Q_{p}(T)V$ for every $p\in \mathbb{H}.$ Now assume that $q \in \Lambda^{S}_{\varepsilon}(B)\setminus \sigma_{S}(B)$. Then 
\begin{equation*}
    \frac{1}{\varepsilon} \leq \|Q_{q}(B)^{-1}\| \leq \|V\| \|V^{-1}\| \|Q_q(T)^{-1}\|.
\end{equation*}
This implies that $\Lambda^{S}_{\varepsilon}(B) \subseteq \Lambda^{S}_{k\varepsilon}(T)$ where $k=\|V\| \|V^{-1}\|.$ The other inclusion follows similarly.

\noindent Proof of (6). Note that $\sigma_{S}(T_{1}\oplus T_{2}) = \sigma_{S}(T_{1}) \bigcup \sigma_{S}(T_{2})$. Further, for every $q \in \mathbb{H}\setminus \sigma_{S}(T_{1}\oplus T_{2})$ we have 
\begin{equation*}
    \|Q_{q}(T_{1}\oplus T_{2})^{-1}\| = \max \big\{\|Q_{q}(T_{1})^{-1}\|,\; \| Q_{q}(T_{2})^{-1} \| \big\}.
\end{equation*}
It follows that $\|Q_{q}(T_{1}\oplus T_{2})^{-1}\| \geq \frac{1}{\varepsilon}$ if and only if either $\|Q_{q}(T_{1})^{-1}\| \geq \frac{1}{\varepsilon}$ or $\|Q_{q}(T_{2})^{-1}\| \geq \frac{1}{\varepsilon}$. Thus the result follows.

\noindent Proof of (7). Clearly $\Lambda^{S}_{\varepsilon}(T)$ is  non-empty as it contains $\sigma_{S}(T)$. Further, it is  closed by continuity of the resolvent map $q \mapsto Q_{q}(T)^{-1}$. Now we show that  $\Lambda^{S}_{\varepsilon}(T)$ is bounded. By the left Cauchy kernel operator series or S-resolvent operator series (see \cite[Definition 4.8.3]{Colombo1}) we have
\begin{equation}\label{Equation: operatorseries}
    \sum\limits_{n\geq 0} T^{n}\; q^{-(n+1)} = -Q_{q}(T)^{-1} (T-\overline{q}I)\; \text{for}\; \|T\|<|q|.
\end{equation}
Further, the operator $(T-\overline{q}I)$ is invertible whenever $\|T\|< |q|$ and the inverse is given by
\begin{equation}\label{Equation: series}
    (T-\overline{q}I)^{-1} = -\sum\limits_{n \geq 0} \big[(\overline{q})^{-1}T\big]^{n} (\overline{q})^{-1} I,
\end{equation}
the series converges in the operator norm. For $\|T\|<|q|$, It follows from Equations  \eqref{Equation: operatorseries} and  \eqref{Equation: series}  that
\begin{align*}
    \|Q_{q}(T)^{-1}\|
    &\leq \Big\| \sum\limits_{n\geq 0} T^{n}\; q^{-(n+1)}\Big\| \; \Big\| \sum\limits_{n \geq 0} \big[(\overline{q})^{-1}T\big]^{n} (\overline{q})^{-1} I \Big\|\\
    &\leq \sum\limits_{n \geq 0} \|T^{n}\| | q^{-(n+1)}| \; \sum\limits_{n\geq 0} \big\|\big[(\overline{q})^{-1}T\big]^{n}\big\| \left|(\overline{q})^{-1}\right|\\
    &\leq\Big[ ~ \sum\limits_{n\geq 0} \frac{\|T\|^{n}}{|q|^{n+1}}~\Big]^{2}\\
    &= \Big[~  \frac{1/|q|}{1-\|T\|/|q|} ~\Big]^{2}\\
    &= \frac{1}{\big[|q|-\|T\| \big]^{2}}\;,
\end{align*}
 Now assume that $|q| > \|T\|+\sqrt{\varepsilon}$. Then $Q_{q}(T)$ is invertible and 
\begin{equation*}
    \|Q_{q}(T)^{-1}\| \leq \frac{1}{\big[|q|-\|T\| \big]^{2}} < \frac{1}{\varepsilon}.
\end{equation*}
Thus $q \notin \Lambda^{S}_{\varepsilon}(T).$ 
Hence 
\begin{align*}
    \Lambda^{S}_{\varepsilon}(T)\subseteq \{q\in \mathbb{H}: |q|\leq \|T\|+\sqrt{\varepsilon}\}.
\end{align*}
In particular,  $\Lambda^{S}_{\varepsilon}(T)$ is bounded.  Thus $\Lambda_{\varepsilon}^{S}(T)$ is a compact subset of $\mathbb{H}$. 
\end{proof}
\section{ Pseudo $S$-spectra of special operators}
In this section, we compute the pseudo $S$-spectra of some special quaternionic operators. To start with, we calculate the pseudo $S$-spectrum of an upper triangular matrix of quaternions. The computation is standard (for instance, see \cite{Cui1} for the complex case). 

\begin{proposition}\label{Prop: upper triangular}
Suppose $\varepsilon>0$ and $A=
\begin{bmatrix}
\lambda_1 & z \\
0 & \lambda_2 
\end{bmatrix}
$, where $\lambda_1, \lambda_2,z \in \mathbb{C}.$ Then
\begin{align*}
    \Lambda^{S}_{\varepsilon}(A)&=\Big \{ q\in \mathbb{H}: \sqrt{\left(|Q_q (\lambda_1)|+|Q_q (\lambda_2)|\right)^2 +\left|(\lambda_1+\lambda_2)z-2\text{Re}(q) z\right |^2+2 |Q_q (\lambda_1)| |Q_q(\lambda_2)|}\\
        & - \sqrt{\left(|Q_q (\lambda_1)|+|Q_q (\lambda_2)|\right)^2 +\left|(\lambda_1+\lambda_2)z-2\text{Re}(q) z\right |^2-2|Q_q (\lambda_1)| |Q_q(\lambda_2)|} \leq 2\varepsilon  \Big \}.
\end{align*}
\end{proposition}
\begin{proof}
For $\mu=\mu_0+i \mu_1 \in \mathbb{C}$, we have
\begin{align*}
Q_\mu(A)&=(A-\mu_0)^2+ \mu_1 ^2\\
&=
\begin{bmatrix}
\lambda_1^2-2\mu_0 \lambda_1 +|\mu|^2 & (\lambda_1+\lambda_2)z-2\mu_0z \\
0 & \lambda_2^2-2\mu_0\lambda_2+|\mu|^2 
\end{bmatrix}.
\end{align*}
Suppose $s_1$ and $s_2$ are singular values of $Q_\mu(A)$ with $s_1\geq s_2$. That is,  $s^{2}_{1}$ and $s^{2}_{2}$ are eigenvalues of $Q_\mu(A)^* Q_\mu(A)$. So, consider the following matrix,   
\begin{align*}
Q_\mu(A)^* Q_\mu(A)=\begin{bmatrix}
|Q_\mu (\lambda_1)|^2 & \overline{Q_{\mu}(\lambda_{1})}[(\lambda_1+\lambda_2)z-2\mu_0z] \\
&\\
{[(\overline{\lambda}_1+\overline{\lambda}_2)\overline{z}-2\mu_{0}\overline{z}]}Q_{\mu}(\lambda_{1}) & |(\lambda_1+\lambda_2)z-2\mu_0 z|^2+|Q_\mu (\lambda_2)|^2 \\
\end{bmatrix}.
\end{align*}
Then
\begin{align*}
s_1^{2}+s_2^{2} &= \text{tr}\left(Q_\mu(A)^* Q_\mu(A)\right)\\
&=|Q_\mu (\lambda_1)|^2+|Q_\mu (\lambda_2)|^2 +|(\lambda_1+\lambda_2)z-2\mu_0 z|^2.
\end{align*}
and 
$s_{1}\; s_{2} = |\det (Q_{\mu}(A))|= |Q_{\mu}(\lambda_{1}) Q_{\mu}(\lambda_{2})|$. Therefore, 
\begin{align*}
(s_{1}\pm s_{2})^2 &= |Q_\mu (\lambda_1)|^2+|Q_\mu (\lambda_2)|^2 +|(\lambda_1+\lambda_2)z-2\mu_0 z|^2 \pm 2 |Q_\mu (\lambda_1)| |Q_\mu(\lambda_2)|. 
\end{align*}
It implies that
\begin{align*}
    s_2 &= \frac{1}{2}\Big[\sqrt{\left(|Q_\mu (\lambda_1)|+|Q_\mu (\lambda_2)|\right)^2 +\left|(\lambda_1+\lambda_2)z-2\mu_0 z\right |^2+2 |Q_\mu (\lambda_1)| |Q_\mu(\lambda_2)|}\\
        &\; - \sqrt{\left(|Q_\mu (\lambda_1)|+|Q_\mu (\lambda_2)|\right)^2 +\left|(\lambda_1+\lambda_2)z-2\mu_0 z\right |^2-2|Q_\mu (\lambda_1)| |Q_\mu(\lambda_2)|}\Big].
\end{align*}
Hence the result follows by using the fact that $\Lambda_{\varepsilon}^{S}(T)$ is axially symmetric and from Proposition \ref{Proposition: singular}. 
\end{proof}
\begin{remark}
If $A$ is a bounded linear operator on a complex  Hilbert space, then
\begin{align*}
\Lambda_\varepsilon(A)=\underset{\|B\|\leq \varepsilon}{\bigcup } \, \sigma(A+B). 
\end{align*}
However, the quaternionic counterpart is not valid.
In particular, for $T,R\in \mathcal{B}(\mathcal{H})$ the one sided inclusion $ \underset{\|S\|\leq \varepsilon}{\bigcup } \, \sigma_S(T+R)\subseteq \Lambda_\varepsilon^S(T)$ 
fails to hold. We provide an example below.
\end{remark}
\begin{example}
Consider $T= \begin{bmatrix}
0 & 1 \\
0 & 0 
\end{bmatrix}
\in M_2(\mathbb{H})$ and $R=\begin{bmatrix}
1.1 \,i & 0 \\
0 & j \\
\end{bmatrix}
\in M_2(\mathbb{H})$. Here $\sigma_{S}(T) = \{0\}$ and $\sigma_{S}(R) = \mathbb{S} \cup \{1.1 m:\; m \in \mathbb{S}\}$. Also $ \|R\|=1.1 $ and $\sigma_S(T+R)=\sigma_{S}(R)$. By Example \ref{Prop: upper triangular} with $\lambda_1=\lambda_2=0$ and $z=1$, we have
\begin{align*}
\Lambda_\varepsilon(T)=\left \{ q\in \mathbb{H}: |q|\leq \left( \varepsilon(\varepsilon+2|q_0|) \right) ^\frac{1}{4} \right \}. 
\end{align*}
In particular, if we choose $\varepsilon = q=1.1$, then  $q \in \sigma_S(T+R)$, but $q\notin \Lambda_\varepsilon(T).$
\end{example}
\begin{remark}
In the complex case, $\Lambda_\varepsilon(\lambda+A)=\lambda+\Lambda_\varepsilon(A)$ for $\lambda\in \mathbb{C}$ and $\varepsilon>0$. But the corresponding quaternionic counter-part $\Lambda_\varepsilon^S(\lambda+T)=\lambda+\Lambda_\varepsilon^S(T)$ fails to hold. Consider $T= \begin{bmatrix}
0 & 1 \\
0 & 0 
\end{bmatrix}$ in Proposition \ref{Prop: upper triangular}. Choose $\varepsilon=1.1$ and $\lambda=1$. Then it is readily seen by following  Example \ref{Prop: upper triangular} that $t:=2.1\in \Lambda_\varepsilon^S(\lambda+T)$ but $t\notin \lambda+ \Lambda_\varepsilon^S(T)$. It is not surprising since $Q_{q}(T)$ is a second order quaternionic operator.
\end{remark}

Before we proceed further, let us recall that if $A$ is a bounded linear operator on a complex Hilbert space, then 
\begin{align}\label{Eq: distance formula}
    \|(A-\lambda I)^{-1}\|\geq  \dfrac{1}{\text{dist}(\lambda, \sigma(A))} \, \, \, \text{ for all } \lambda \in \mathbb{C} \setminus \sigma(A), 
\end{align}
where $\text{dist} (\lambda,\sigma(A))=\underset{\mu\in \sigma(A)}{\inf}\, |\lambda-\mu|$. Moreover, the equality holds for a certain class of operators like normal, subnormal, hyponormal etc.
It is evident from the above distance formula \eqref{Eq: distance formula} that the $\varepsilon$-pseudospectrum of $A$ contains the closed disc of radius $\varepsilon$ around the spectrum $\sigma(A)$ \cite{Tref2}.

\subsection{A distance formula}
In the quaternionic setup, a natural analogue of the distance formula is available  for quaternionic normal operators. Further, an important estimate of the pseudo resolvent operator is obtained in  \cite{Paula} for normal operators using slice functional calculus. We recall the distance formula here.  

\begin{lemma}\label{lemma: distance formula}
Let $T\in \mathcal{B}(\mathcal{H})$ be normal. Then 
\begin{align}\label{Equation: normcondn}
    \|Q_{q}(T)^{-1}\|=  \frac{1}{\inf \big\{ |(\mu- q) (\mu - \overline{q}) |:\; \mu \in \sigma_{S}(T)\cap \mathbb{C} \big\}},
\end{align}
for every $q \in \rho_{S}(T)$.
\end{lemma}
\begin{proof}
Without loss of generality, we consider $\lambda \in \rrho_{s}(T)\cap \mathbb{C}$. Since $T$ is normal, we know that $Q_{\lambda}(T)$ is normal and so its inverse. Further, we have  
\begin{equation}\label{Equation: distance}
    \|Q_{\lambda}(T)^{-1}\| = r_{s}(Q_{\lambda}(T)^{-1}) = \sup\big\{ |p|:\; p \in \sigma_{S}(Q_{\lambda}(T)^{-1})\big\}.
\end{equation}
Now our task is to express the $\sigma_{S}(Q_{\lambda}(T)^{-1})$ in terms of $\sigma_{S}(T)$. To obtain this, let us define a function $f_{\lambda}(x) =\left( x^2-2 \text{Re}(\lambda)x+|\lambda|^2 \right)^{-1},\; \text{for all}\; x \in \mathbb{H}$ with $|x|\leq \|T\|.$ Clearly, $f_{\lambda}(x) \in \mathbb{C}_{m}$ whenever $x \in \mathbb{C}_{m}$, $m \in \mathbb{S}$. Further, $f_{\lambda}(x)$ is a $S$-regular (slice) function since for every $m \in \mathbb{S}$ and $\alpha, \beta \in \mathbb{R}$ we have
\begin{equation*}
    \frac{\partial}{\partial \alpha} f_{\lambda}(\alpha + m \beta) + m \frac{\partial}{\partial \beta} f_{\lambda}(\alpha + m \beta) = 0\; \;\text{and}\;\;  \frac{\partial}{\partial \alpha} f_{\lambda}(\alpha + m \beta) +  \frac{\partial}{\partial \beta} f_{\lambda}(\alpha + m \beta) m = 0.
\end{equation*}
It implies that $f$ is $\mathbb{H}$-intrinsic slice function. 
From the spectral mapping theorem for $\mathbb{H}$-intrinsic slice functions
we see that
\begin{equation*}
    \sigma_{S}(Q_{\lambda}(T)^{-1}) = \sigma_{S}(f_{\lambda}(T)) = \Big\{ f_{\lambda}(s):\; s \in \sigma_{S}(T)\Big\}.
\end{equation*}
Therefore, from Equation \eqref{Equation: distance} it follows that 
\begin{align*}
    \|Q_{\lambda}(T)^{-1}\| &= \sup\Big\{ |f_{\lambda}(s)|:\; s \in \sigma_{S}(T)\Big\}\\
    &=  \sup\Big\{ \frac{1}{|s^{2} - 2 Re(\lambda)s+|\lambda|^{2}|}:\; s \in \sigma_{S}(T)\Big\}\\
    &= \sup\Big\{ \frac{1}{|\mu^{2} - 2 Re(\lambda)\mu+|\lambda|^{2}|}:\; \mu \in \sigma_{S}(T)\cap \mathbb{C}\Big\}\\
    &= \frac{1}{\inf \big\{ |(\mu-\lambda) (\mu-\overline{\lambda})|:\; \mu \in \sigma_{S}(T)\cap \mathbb{C} \big\}}.
\end{align*}
Finally, from the definition of $Q_{q}(T)$, the result follows for every $q\in \rrho_{S}(T).$
\end{proof}

 In the current scenario, we introduce the notion of $\varepsilon$-sphere around the $S$-spectrum. Later we shall show that for certain class of quaternionic operators the pseudo $S$-spectrum coincides with $\varepsilon$-sphere.

\begin{definition}\label{Definition: epsilonsphere}
For $\varepsilon > 0$, the $\varepsilon$-sphere around a compact set $K$ in $\mathbb{H}$  is defined by
\begin{equation*}
    S^{\mathbb{H}}(K, \varepsilon) = \big\{q\in \mathbb{H}: \underset{\mu \in K\cap \mathbb{C}}{\inf} \,\, \big|\mu^2- 2 \text{Re}(q)\mu+|q|^2 \big| \leq \varepsilon \big\}.
\end{equation*}
Notice that $S^{\mathbb{H}}(K, \varepsilon)$ is axially symmetric and $\bigcap\limits_{\varepsilon>0} S^{\mathbb{H}}(K, \varepsilon) = K.$ 
\end{definition}
The case when $K=\sigma_S(T)$ is of special interest for understanding the structure of pseudo $S$-spectra for a large class of quaternionic operators $T\in \mathcal{B}(\mathcal{H})$.

Now we turn our discussion to point out some of the important properties of $S^{\mathbb{H}}(\sigma_{S}(T), \varepsilon)$ and obtain its relation with the pseudo $S$-spectrum. To begin with, $S^{\mathbb{H}}(\sigma_{S}(T), \varepsilon)$ contains $\sigma_{S}(T)$ but may not intersect the real line for some $\varepsilon$. See the example below. 
\begin{example}
Let 
$A = \begin{bmatrix}{}
j & 0\\ 0 & k
\end{bmatrix}\in M_{2}(\mathbb{H})$.
Then $\sigma_{S}(A)  = \mathbb{S}$. A real number $r \in S^{\mathbb{H}}(\mathbb{S},\; \varepsilon) \cap \mathbb{R}$ if and only if 
\begin{align*}
    \varepsilon \geq \inf\limits_{m \in \mathbb{S}} |m^{2}-2rm+r^{2}| = \inf\limits_{m \in \mathbb{S}} |m-r|^{2} =\inf\limits_{m\in \mathbb{S}} (m-r)(\overline{m}-r) =  \inf\limits_{m \in \mathbb{S}} (|m|^{2}+r^{2})= 1+r^2 .
\end{align*}
This inequality is valid if and only if $\varepsilon \geq 1$. Hence  $ S^{\mathbb{H}}(\mathbb{S},\; \varepsilon) \cap \mathbb{R} = \emptyset$ whenever $\varepsilon<1.$
\end{example}
In the following theorem, we show that, like in the complex case, for a quaternionic normal operator $T \in \mathcal{B}(\mathcal{H})$ the pseudo $S$-spectrum of $T$ coincides with $\varepsilon$-sphere around $\sigma_{S}(T)$ for every $\varepsilon>0.$ 
\begin{theorem}\label{normal thm}
Let $T\in \mathcal{B}(\mathcal{H})$ be normal and $\varepsilon>0$.
Then 
$$S^\mathbb{H}(\sigma_S(T),\varepsilon)= \Lambda^S_\varepsilon (T).$$ 

\end{theorem}

\begin{proof} 
Suppose that $T$ is a normal operator. Then $Q_{q}(T)^{-1}$ is normal and so $\|Q_{q}(T)^{-1}\|=r_S(Q_{q}(T)^{-1}).$ Consequently, the equality follows from Lemma \ref{lemma: distance formula}.
\end{proof}

Since $\sigma_{S}(T)$ is axially symmetric, $\sigma_{S}(T)\cap \mathbb{C}$ is a non-empty set and each equivalence class in the $S$-spectrum is uniquely determined by the members of  $\sigma_{S}(T)\cap \mathbb{C}$.
The following theorem gives a relationship between the $\varepsilon$-sphere around $\sigma_{S}(T)$ and certain closed disks around the members of $\sigma_{S}(T)\cap \mathbb{C}$.

\begin{lemma} \label{Lemma: unionofdisks}
Let $T \in \mathcal{B}(\mathcal{H})$ and $\varepsilon > 0$. Then 
\begin{equation*}
    S^{\mathbb{H}}(\sigma_{S}(T), \varepsilon) \subseteq \bigcup\limits_{\lambda\in \sigma_{S}(T)\cap \mathbb{C}}\Omega_{D\big(\text{Re}(\lambda),\; \sqrt{\varepsilon+d^{2}}\big)},
\end{equation*}
where $d = dist(0, \text{Im}(\sigma_{S}(T)\cap \mathbb{C})).$
\end{lemma}
\begin{proof}
Let $q \in S^\mathbb{H}(\sigma_S(T),\varepsilon) $. Then $\inf\limits_{\lambda \in \sigma_{S}(T)\cap \mathbb{C}} \big|\lambda^{2}- 2 \text{Re}(q)\lambda+|q|^{2} \big|\leq \varepsilon$. If we take $\mu = \text{Re}(q)+i |\text{Im}(q)|$, it follows that 
\begin{align*}
    \varepsilon & \geq  
     \inf \limits_{\lambda \in \sigma_{S}(T)\cap \mathbb{C}} \big|\lambda^{2}- 2 \text{Re}(\mu)\lambda+|\mu|^{2} \big| \\
    &= \inf \limits_{\lambda \in \sigma_{S}(T)\cap \mathbb{C}} \big|(\mu-\lambda ) (\mu - \overline{\lambda}) \big|\\
    &= \inf \limits_{\lambda \in \sigma_{S}(T)\cap \mathbb{C}} \big|\big(\mu-\text{Re}(\lambda)\big)^{2}+ |\text{Im}(\lambda)|^{2} \big|.
\end{align*}
Since $\sigma_{S}(T)$ is compact, there exists a $\lambda_{0} \in \sigma_{S}(T)\cap \mathbb{C}$ such that   
\begin{align*}
    |\mu - \text{Re}(\lambda_{0})|^{2} &= \inf\limits_{\lambda \in \sigma_{S}(T)\cap \mathbb{C}}|\mu-\text{Re}(\lambda)|^{2}\\
    &= \inf\limits_{\lambda \in \sigma_{S}(T)\cap \mathbb{C}}\big|\big(\mu-\text{Re}(\lambda)\big)^{2} + |\text{Im}(\lambda)|^{2} - |\text{Im}(\lambda)|^{2}\big|\\
    &\leq \inf\limits_{\lambda \in \sigma_{S}(T)\cap \mathbb{C}}\big|\big(\mu-\text{Re}(\lambda)\big)^{2} + |\text{Im}(\lambda)|^{2} \big| + \inf\limits_{\lambda \in \sigma_{S}(T)\cap \mathbb{C}}|\text{Im}(\lambda)|^{2}\\
    &\leq \varepsilon+  d^{2},
\end{align*}
where $d= dist(0, \text{Im}(\sigma_{S}(T)\cap \mathbb{C})).$ Hence $\mu \in \bigcup\limits_{\lambda \in \sigma_{S}(T)\cap \mathbb{C}} D\big( \text{Re}(\lambda),\; \sqrt{\varepsilon+d^{2}}\big)$ and the desired inclusion follows.  
\end{proof} 
\noindent Now we give a concrete description of pseudo $S$-spectrum of a self-adjoint operator. 
\begin{theorem}\label{Self adjoint S-spectrum}
Let $T\in \mathcal{B}(\clh)$ be self-adjoint and $\varepsilon>0$. Then
\begin{align*}
    \Lambda_\varepsilon^S(T)=\bigcup\limits_{\lambda \in \sigma_{S}(T)}\Omega_{D(\lambda,\;\sqrt{\varepsilon})}.
\end{align*}

\end{theorem}
\begin{proof} Since $T$ is self-adjoint, $\sigma_{S}(T) \subset \mathbb{R}$ and so $d = 0.$ By Lemma \ref{Lemma: unionofdisks} and Theorem \ref{normal thm} we see that 
\begin{equation*}
    \Lambda_\varepsilon^S(T)= S^{\mathbb{H}}(\sigma_{S}(T), \varepsilon) \subseteq \bigcup\limits_{\lambda \in \sigma_{S}(T)}\Omega_{D(\lambda,\;\sqrt{\varepsilon})}.
\end{equation*}
Now we prove the reverse inclusion. Let  $q \in \Omega_{D(r,\sqrt{\varepsilon})}$ for some $r \in \sigma_{S}(T)$. If we take $z= \text{Re}(q) + i|\text{Im}(q)|$ then $ z \in D(r, \sqrt{\varepsilon})$. Further, 
\begin{equation*}
    \inf\limits_{t \in \sigma_{S}(T)}|z^{2}-2t z+t^{2}| \leq |z-r|^{2}\leq \varepsilon,
\end{equation*}
 which yields  $z \in S^{\mathbb{H}}(\sigma_S(T),\varepsilon)$. The reverse inclusion follows from the fact that $S^{\mathbb{H}}(\sigma_{S}(T), {\varepsilon})$ is axially symmetric.
 \end{proof}
In the following theorem, we characterize orthogonal projections through their pseudo $S$-spectra. 
\begin{theorem}\label{Thm: projection}
Let $T\in \mathcal{B}(\mathcal{H})$ be a (non-trivial) normal operator. Then 
\begin{align*}
T \text{ is an orthogonal projection if and only if 
} \Lambda_\varepsilon^S(T)= \Omega_{D(0,\sqrt{\varepsilon})} \cup \Omega_{D(1,\sqrt{\varepsilon})} \, \, \text{ for all } \varepsilon>0.
\end{align*}

\end{theorem}
\begin{proof}
Suppose $T$ is a non-trivial orthogonal projection. We note that $ \{0,1 \} \subseteq \sigma_S(T)$. Indeed, $Q_0(T)=T$, which is not invertible since $T\neq I$. Also, $Q_1(T)=I-T$, which is not invertible since $T\neq 0$. 

Let $q\in \mathbb{H} \setminus \{0,1\}$. Then $Q_{q}(T)=(1-2\text{Re}(q))T+|q|^2I$. Take $\lambda=1-2\text{Re}(q)$. An elementary calculation shows that 
\begin{align*}
Q_{q}(T)^{-1}=\frac{1}{|q|^2}\left(I-\frac{\lambda T}{|q|^2+\lambda}\right).
\end{align*}
Here note that $|q|^2+\lambda \neq 0$ since $q\neq 1$. This gives 
$$\sigma_S(T)=\{0,1 \}.$$
Since T is self adjoint, by Theorem \ref{Self adjoint S-spectrum}, it follows that 
\begin{align*}
\Lambda_{\varepsilon}^S(T)=\Omega_{D(0,\sqrt{\varepsilon})} \cup \Omega_{D(1,\sqrt{\varepsilon})} \, \, \text{ for all } \varepsilon>0.
\end{align*}

Conversely, let us assume that $\Lambda_{\varepsilon}^S(T)=\Omega_{D(0,\sqrt{\varepsilon})} \cup \Omega_{D(1,\sqrt{\varepsilon})} \, \, \text{ for all } \varepsilon>0.$
Then by Proposition \ref{Elementary prop}(3), we have
\begin{align*}
\sigma_S(T)=\{0,1\}.
\end{align*}
Since $T$ is normal and $\sigma_S(T)\subseteq \mathbb{R}$, $T$ is self adjoint.
By polynomial functional calculus for quaternions, it follows that
\begin{align*}
\sigma_S(T-T^2)=\big\{[\lambda-\lambda^2]:\; \lambda\in \sigma_S(T)\cap \mathbb{C}\big\}=\{0\}.
\end{align*}
Thus $r_S(T-T^2)=0$. Since $T$ is self adjoint, we have $\|T-T^2\|=r_S(T-T^2)=0$. Hence $T^2=T$ and so $T$ is an orthogonal projection.
\end{proof}

\subsection{\bf Left multiplication operator} Now we turn our discussion to left multiplication in right quaternionic Hilbert space. To begin with, let $\mathcal{H}$ be a separable right quaternionic Hilbert space with an orthonormal basis $\{e_{n}:\; n\in \mathbb{N}\}$.
For $q \in \mathbb{H}$, the \textit{left multiplication operator} $L_{q}\colon \mathcal{H} \to \mathcal{H}$ \textit{induced by $q$} is defined as
\begin{equation*}
    L_{q}(x) = q \cdot x = \sum\limits_{n\in \mathbb{N}} e_{n} \cdot q \langle e_{n}, x\rangle \; \,\,\,\, (x\in \mathcal{H}).
\end{equation*}

The norm of $L_{q}$ is computed as follows:
\begin{equation*}
    \|L_{q}\| = {\|L_{q}^{\ast}L_{q}\|}^{1/2} = \|L_{|q|^{2}}\|^{1/2} = \Big[\sup\limits_{\|x\|=1}  \sum\limits_{n=1}^{\infty} e_{n}|q|^{2} \langle e_{n},\; x\rangle\Big]^{1/2} = |q|.
\end{equation*} 

\begin{proposition} \label{Prop: Mult op real}
Let $q\in \mathbb{H}$ and $\varepsilon>0$. Then $\Lambda_{\varepsilon}^{S}(L_{q}) = S^{\mathbb{H}}([q], \varepsilon)$. In particular, if $q\in \mathbb{R}$, then $\Lambda_{\varepsilon}^{S}(L_{q})= \Omega_{D(q,\sqrt{\varepsilon})}$.
\end{proposition}
\begin{proof}
The adjoint of the operator $L_{q}$  is $L_q^\ast=L_{\overline{q}}$ and so $L_{q}$ is a normal operator. We now compute the $S$-spectrum of $L_{q}.$ Note that
\begin{align*}
    p \in \sigma_{s}(L_{q}) &\iff Q_{p}(L_{q})\; \text{is not invertible in}\; \mathcal{B}(\mathcal{H})\\
    &\iff L_{q^{2}} - 2 Re(p)L_{q}+|p|^{2}I \; \text{is not invertible in}\; \mathcal{B}(\mathcal{H})\\
    &\iff L_{(q^{2}-2 Re(p)q+|p|^{2})} \; \text{is not invertible in}\; \mathcal{B}(\mathcal{H})\\
    &\iff q^{2}-2 Re(p)q+|p|^{2} = 0\\
    & \iff p \in [q].
\end{align*}
Therefore, $\sigma_{S}(L_{q}) = [q].$ Then by Theorem \ref{normal thm}, $\Lambda_{\varepsilon}^{S}(L_{q}) = S^{\mathbb{H}}([q], \varepsilon)$. Finally,  if $q$ is a real number, then $L_q$ is self adjoint. Then the final result follows directly from Theorem \ref{Self adjoint S-spectrum}.
\end{proof}


In the following theorem we discuss the converse of the above proposition. 
\begin{theorem}\label{bound for s-pseudo}
Let $T \in \mathcal{B}(\mathcal{H})$ and $q \in \mathbb{H}$ be such that $\Lambda^{S}_{\varepsilon}(T) = S^{\mathbb{H}}([q], \varepsilon)$ for every $\varepsilon >0$. Then the following holds:
\begin{enumerate}
    \item  If $q = r\in \mathbb{R}$ then $T = L_{r}.$
    \item If $q \in \mathbb{H}\setminus \mathbb{R}$ and $T$ is normal then there exists a unitary operator $U$ on $\mathcal{H}$ satisfying,
    \begin{equation*}
        T = U^{\ast}L_{q}U.
    \end{equation*}
\end{enumerate}
\end{theorem}

\begin{proof}
 Proof of $(1).$ If $q=r \in \mathbb{R}$ then  $\Lambda^{S}_{\varepsilon}(T) = S^{\mathbb{H}}(r, \varepsilon)= \Omega_{D(r, \sqrt{\varepsilon})}$ for every $\varepsilon >0.$ If $\mathcal{W} = \Omega_{D(r, {\varepsilon})}$ then $\mathcal{W}$ is an axially symmetric $s$-domain containing $\sigma_{S}(T)$ and hence it is a $T$-admissable open set. It is evident \allowdisplaybreaks[1] that $\partial(\mathcal{W}\cap \mathbb{C}) = \big\{z\in \mathbb{C}:\; |z-r| = {{\varepsilon}}\big\} $. By the quaternionic functional calculus \cite{Colombo1}, we have 
\begin{align*}
    \|T-L_{r}\| &= \frac{1}{2\pi }\Big\| \int\limits_{\partial(\mathcal{W}\cap \mathbb{C})} S_{L}^{-1}(s, T)\; \frac{ds}{i} (s-r)\Big\|\\
    &\leq \frac{1}{2\pi}\; \int\limits_{0}^{2\pi} \|Q_{s}(T)^{-1}\| \|(T-\overline{s}I)\||s-r|\; ds\\
    &\leq \frac{1}{2\pi}\; \int\limits_{0}^{2\pi} \|Q_{s}(T)^{-1}\| (\|(T\|+|s|)|s-r|\; ds\\
    &< \frac{1}{2\pi {\varepsilon}} \varepsilon \int\limits_{0}^{2\pi} (\|T\|+\varepsilon+r)\; ds\\
    &=\frac{1}{2\pi} (\|T\|+\varepsilon+r) 2\pi \varepsilon\\
    &=(\|T\|+\varepsilon+r)\varepsilon.
\end{align*}
Since $\varepsilon>0$ is arbitrary, it follows that $T = L_{r}.$

\noindent Proof of (2). Since $T$ is normal, by Theorem 9.3.5 and Remark 9.3.7 of \cite{Colombo2}, there exists an anti self-adjoint unitary operator $J \in \mathcal{B}(\mathcal{H})$ that commutes with $T$ such that $T = A + J B$ where $A$ is a self-adjoint operator, $B$ is a positive operators. Originally, this decomposition of normal operator is due to O. Teichm\"{u}ller \cite{Teichmuller}. Define $\mathcal{H}^{Ji}_{+}:= \big\{x\in \mathcal{H}:\; J(x) = xi \big\}$. Then $\mathcal{H}^{Ji}_{+}$ is a complex Hilbert space with the inner product induced from $\mathcal{H}$ and it is invariant under $T$. 
Let $T_{+}:= T\big|_{\mathcal{H}^{Ji}_{+}}$. Then it is a bounded complex linear on $\mathcal{H}^{Ji}_{+}$ such that 
\begin{equation*}
    T(x) = T_{+}(x_{+}) - T_{+}(x_{-}\cdot j)\cdot j, \; \text{for every}\; x = x_{+}+x_{-} \in \mathcal{H}^{Ji}_{\pm}.
\end{equation*}
By 
(3) of Proposition \ref{Elementary prop}, it follows that
\begin{align*}
    \sigma(T_{+}) = \sigma_{S}(T) \cap \mathbb{C}^{+} = \bigcap\limits_{\varepsilon >0}\Lambda_{\varepsilon}^{S}(T) \cap \mathbb{C}^{+} = \bigcap\limits_{\varepsilon >0} S^{\mathbb{H}}([q], \varepsilon)\cap \mathbb{C}^{+} = \big\{z:= \text{Re}(q)+i|\text{Im}(q)| \big\}.
\end{align*}
Thus  $\Lambda_{\varepsilon}(T_{+}) = D(z, \varepsilon)$ for every $\varepsilon>0.$ Further, by \cite[Corollary 3.17]{Krishnan} it implies that $T_{+} = z\cdot I_{\mathcal{H}^{Ji}_{+}}$. Now suppose that $\{f_{n}:\; n\in \mathbb{N}\} \subset \mathcal{H}^{Ji}_{+}$ is a Hilbert basis for $\mathcal{H}^{Ji}_{+}$. Then it is also a Hilbert basis for $\mathcal{H}$ and hence for every $x = x_{+} + x_{-}$ where $x_{\pm} \in \mathcal{H}^{Ji}_{\pm}$ we have 
\begin{align*}
    T(x) &= T_{+}(x_{+}) - T_{+}(x_{-}\cdot j) \cdot j\\
    &= z \cdot x_{+} - z \cdot  (x_{-}\cdot j) \cdot j\\
    & = \sum\limits_{n \in \mathbb{N}} f_{n}\cdot z \langle f_{n},\; x_{+}\rangle - \sum\limits_{n \in \mathbb{N}} f_{n}\cdot z \langle f_{n},\; x_{-} \cdot j\rangle \cdot j\\
    & = \sum\limits_{n \in \mathbb{N}} f_{n}\cdot z \langle f_{n},\; x\rangle\\
    & = L_{z}(x).
\end{align*}
Note that  $s: = \frac{\sqrt{q_{1}^{2}+q_{2}^{2}+q_{3}^{2}}+q_{1}-q_{3} j + q_{2}k}{|\sqrt{q_{1}^{2}+q_{2}^{2}+q_{3}^{2}}+q_{1}-q_{3} j + q_{2}k|}$ is a unit quaternion with $\overline{s}qs = z$. Thus the operator defined by $U = L_{s}$ is unitary. By direct computation, we get  $U^{\ast}L_{q}U = T.$
\end{proof}
\begin{theorem} \label{Theorem: emptyorconnected}
Let $q \in \mathbb{H}$. Then $\Lambda^{S}_{\varepsilon}(L_{q})\cap \mathbb{R}$ is an empty set if $\varepsilon< |\text{Im}(q)|^{2}$. Otherwise, $\Lambda^{S}_{\varepsilon}(L_{q})\cap \mathbb{R}$  is a non-empty connected set. 
\end{theorem}
\begin{proof}
From the definition of $\Lambda^{S}_{\varepsilon}(L_{q})$, it is immediate to see that a real number $r \in \Lambda^{S}_{\varepsilon}(L_{q})$ if and only if 
\begin{equation} \label{Equation: inequality}
    \big| r- [\text{Re}(q)+i |\text{Im}(q)|] \big|^{2} = (r-\text{Re}(q))^{2}+|\text{Im}(q)|^{2}= \big| r^{2} - 2 \text{Re}(q)r+|q|^{2}\big| \leq \varepsilon.
\end{equation}
 If $\varepsilon < |\text{Im}(q)|^{2}$ then  by above observation, $(r-\text{Re}(q))^{2}+ |\text{Im}(q)|^{2} < |\text{Im}(q)|^{2}$. It shows that there does not exist any $r\in \mathbb{R}$ satisfying the above inequality. In this case, $\Lambda_{\varepsilon}^{S}(L_{q}) \cap \mathbb{R}$ is an empty set. 
 
 On the other hand if $\varepsilon\geq |\text{Im}(q)|^{2}$, then $\text{Re}(q)$ satisfies the Equation \eqref{Equation: inequality} and so $\Lambda_{\varepsilon}^{S}(L_{q})\cap \mathbb{R}$ is a non-empty set. Let 
 $r_{1}, r_{2} \in \Lambda_{\varepsilon}^{S}(L_{q})\cap \mathbb{R}$. Then  
\begin{equation*}
    |r_{\ell}- (\text{Re}(q)+i |\text{Im}(q)|)|  \leq \sqrt{\varepsilon}\; \; \text{for}\; \ell = 1,2.
\end{equation*}
Finally,  for any $t \in (0, 1)$,  we have
\begin{align*}
    \big|r_{1}t + r_{2}(1-t)-&[\text{Re}(q)+i|\text{Im}(q)|]\big|\\
    &= \big|r_{1}t + r_{2}(1-t)-[\text{Re}(q)+i|\text{Im}(q)|] (t+1-t)\big|\\
    &\leq \big|r_{1}-[\text{Re}(q)+i|\text{Im}(q)|]\big|t + \big|r_{2}-[\text{Re}(q)+i|\text{Im}(q)|]\big|(1-t)\\
    & \leq \sqrt{\varepsilon}t+ \sqrt{\varepsilon}(1-t)\\
    &\leq \sqrt{\varepsilon}. 
\end{align*}
Hence $\Lambda^{s}_{\varepsilon}(L_{q})\cap \mathbb{R}$ is either an empty set or connected whenever it is non-empty.
\end{proof}
The above result may not hold true for an arbitrary  quaternionic operator. We give an example of a normal matrix with quaternion entries for which the intersection of pseudo $S$-spectrum and the real line is non-empty but disconnected.
\begin{example}
Consider the normal matrix  $A = \begin{bmatrix}
(1+\frac{1}{2}j) & 0\\ 0 & 3-\frac{1}{2\sqrt{2}}(i+j)
\end{bmatrix} \in M_{2}(\mathbb{H})$ then $\sigma_{S}(A) = \big\{1+\frac{1}{2}m,\; 3+\frac{1}{2}m:\; m \in \mathbb{S} \big\}.$ If we take  $\varepsilon = \frac{1}{4}$, then the real scalars $1, 3 \in \Lambda_{\frac{1}{4}}^{S}(A)$ but $2 \notin \Lambda_{\frac{1}{4}}^{S}(A)$ because $\min \big\{|2-(1+\frac{1}{2}i)|,\; |2-(3+\frac{1}{2}i)|\big\} = \frac{\sqrt{5}}{2} > \frac{1}{2}$. Hence $\Lambda_{\frac{1}{4}}^{S}(A)\cap \mathbb{R}$ is a nonempty set but disconnected. 
\end{example}
\subsection*{Comments on the structure of $\Lambda^{S}_{\varepsilon}(L_{q})$} Now we discuss the structure of pseudo $S$-spectra of left mulitiplication operators via the following example. Contrary to complex Hilbert spaces, the structure of  $\Lambda^{S}_{\varepsilon}(L_{q})$ differs drastically for various values of $\varepsilon$.  
\begin{example}
If $q = \frac{1}{2}(1+\frac{i}{5})$ then $\Lambda_{\varepsilon}^{S}(L_{q}) = S^{\mathbb{H}}([\frac{1}{2}(1+\frac{i}{5})], \varepsilon)$. So from the Definition \ref{Definition: epsilonsphere} we see that
\begin{equation*}
    \Lambda^{S}_{\varepsilon}(L_{q}) = \Omega_{\big\{\lambda \in \mathbb{C}\colon\; |(\lambda-\frac{1}{2}-\frac{i}{5})(\lambda-\frac{1}{2}+\frac{i}{5})| \leq \varepsilon\big\}}.
\end{equation*}
Equivalently, for a complex number $\lambda := x+iy$, the equivalence class $[\lambda] \in   \Lambda^{S}_{\varepsilon}(L_{q})$ if and only if 
\begin{equation}\label{Equation: SH}
    \Big[ \big(x-\frac{1}{2}\big)^{2} + \big(\frac{1}{100}-y^{2}\big)\Big]^{2} + 4 \big(x-\frac{1}{2}\big)^{2} \leq \varepsilon^{2}.
\end{equation}
Now we plot the graph of Equation \eqref{Equation: SH} using MATLAB for various value of $\varepsilon$. Firstly, consider the value of $\varepsilon$ less than $|\text{Im}(q)|$ i.e., $\varepsilon < 0.1 $ but close to $0.1$ then the pseudo $S$-spectrum of $L_{q}$ has an ``eight" like shape when it is projected to the complex plane as shown in Figure (A).  If we choose  $\varepsilon \geq 0.1$ then it represents an ellipse like structure in Figure (B) and it appears that when $\varepsilon$ is way bigger than $0.1$, the pseudo $S$-spectrum of $L_{q}$ turns out to be a perfect ellipsoid. In the Figure (C), we notice that the eight like shape breaks into a pair of ellipses that are symmetric about $X$-axis when $\varepsilon$ is close to zero.  

\begin{figure}[h]
    \centering
    \subfloat[\centering when $\varepsilon <0.1$]{{\includegraphics[width=7cm]{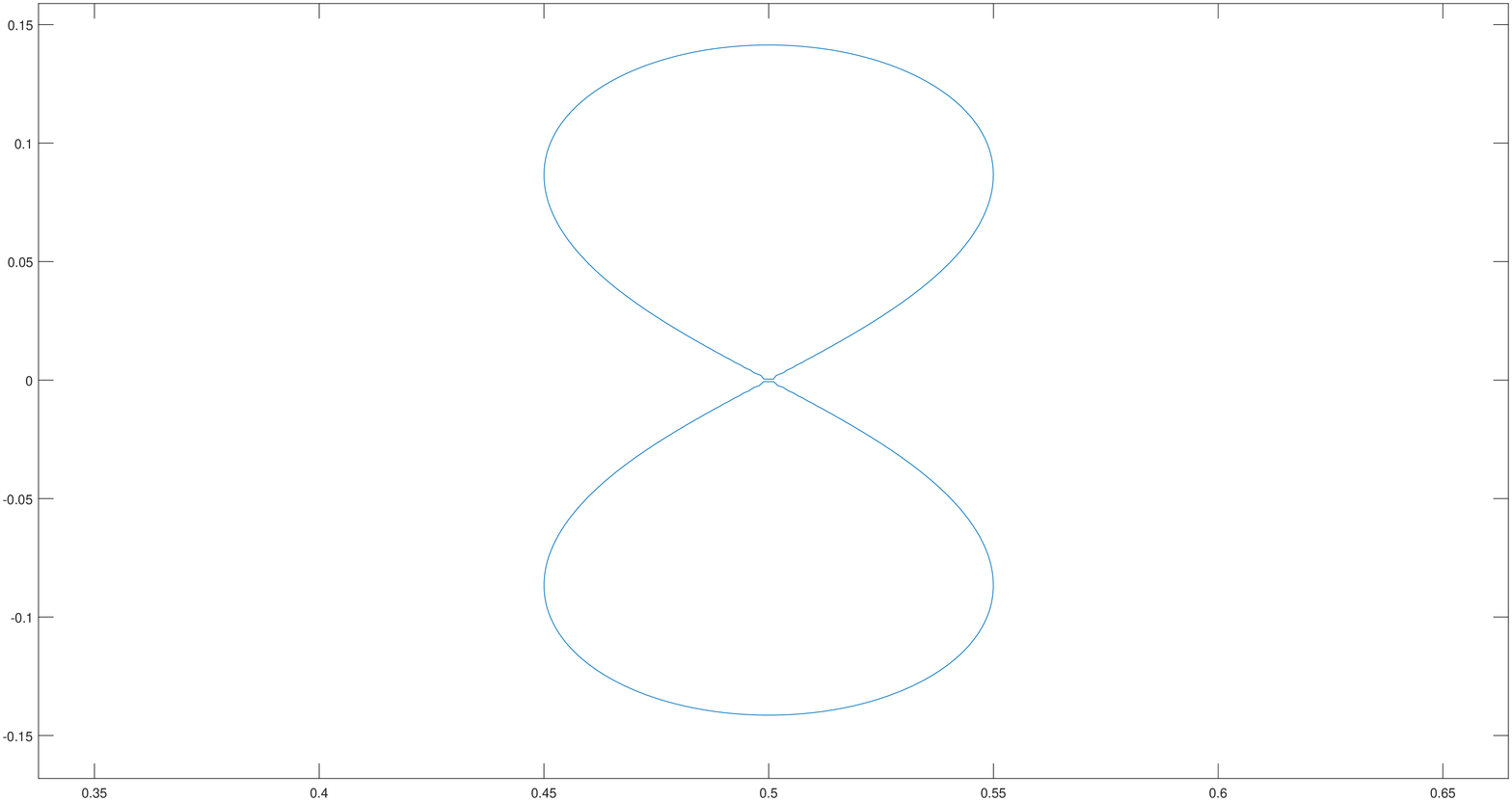} }}%
    \qquad
    \subfloat[\centering when $\varepsilon \geq 0.1$]{{\includegraphics[width=7cm]{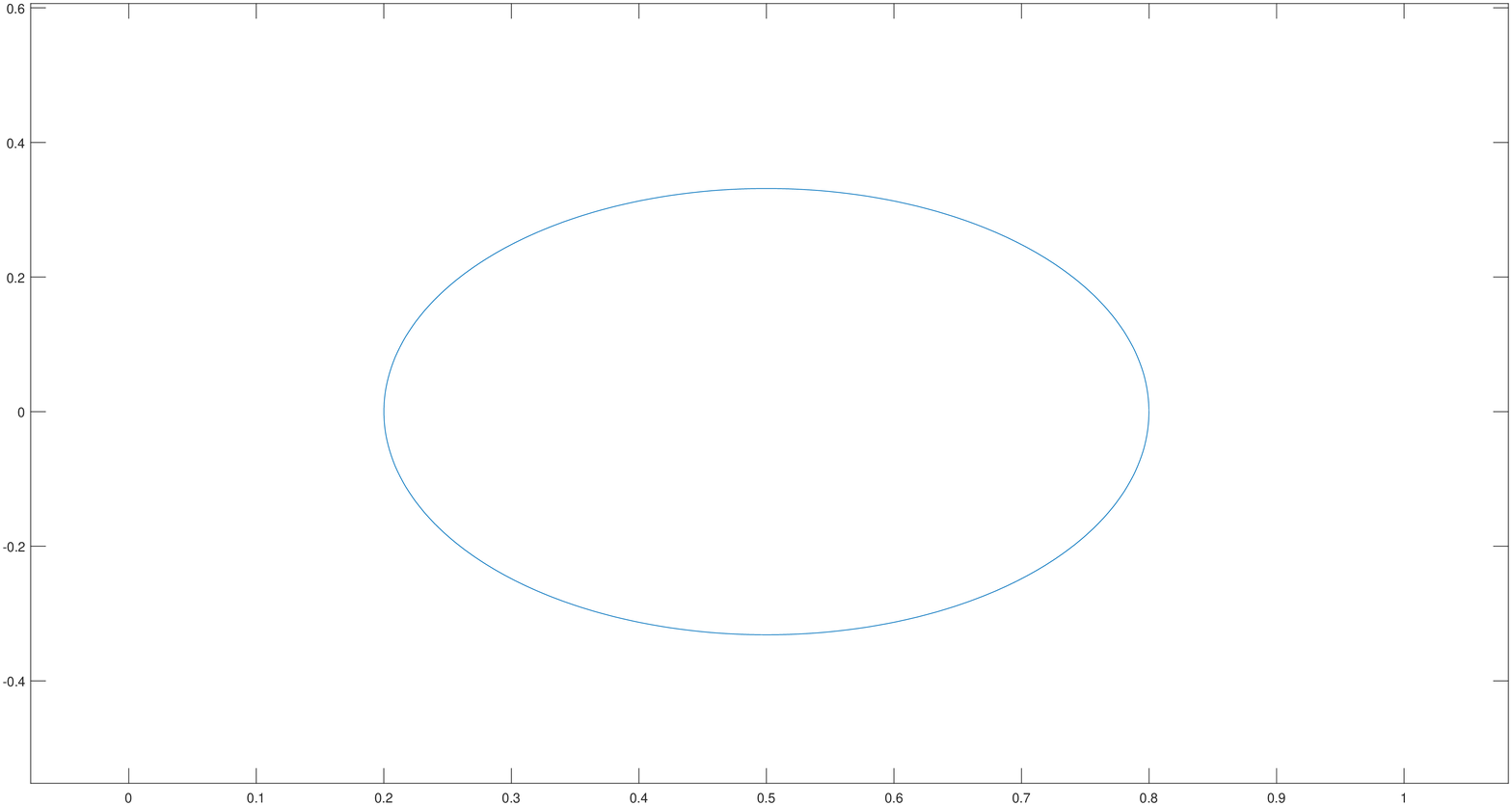} }}%
     \subfloat[\centering when $\varepsilon <<0.1$]{{\includegraphics[width=7cm]{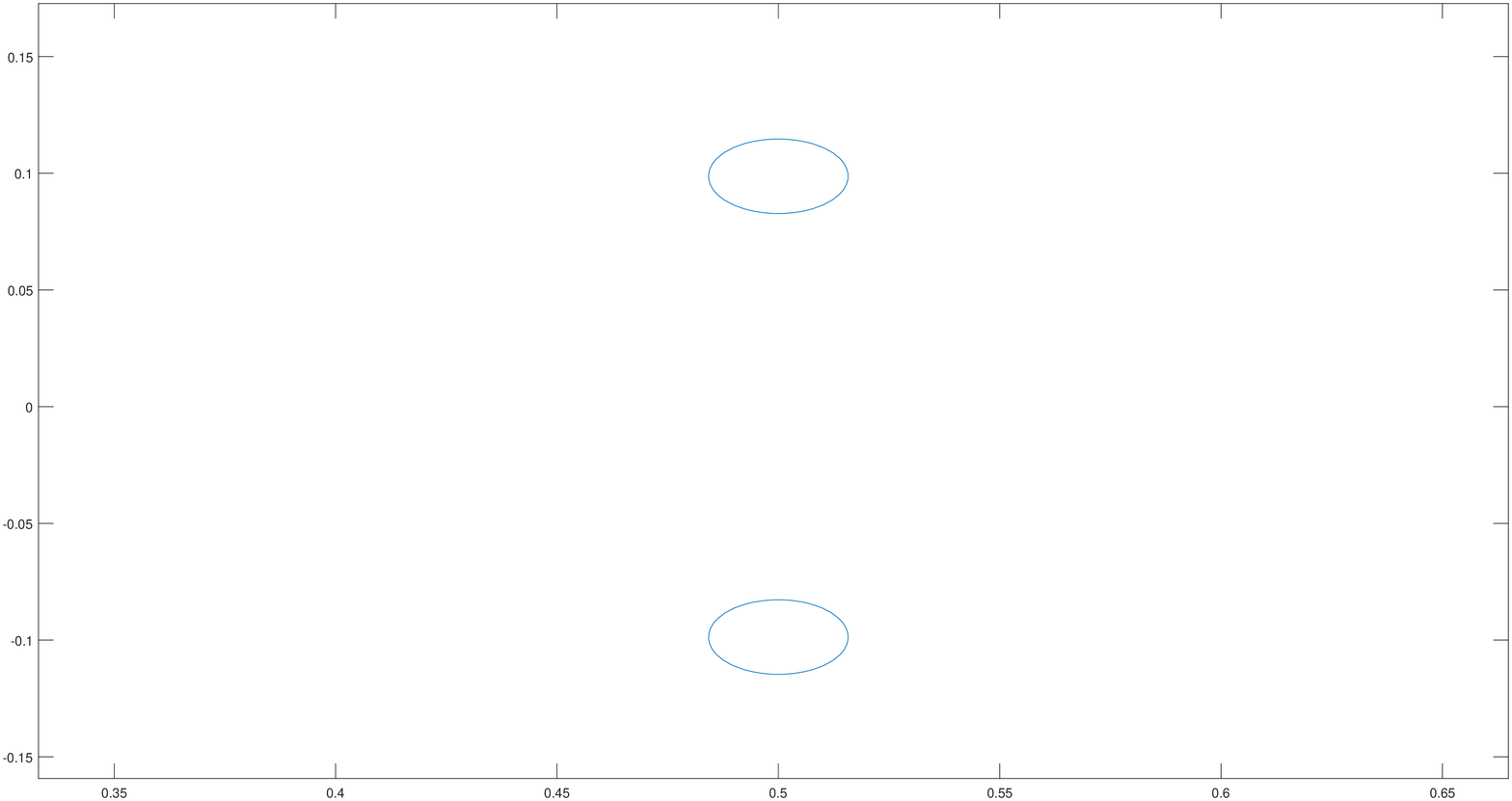} }}%
    \label{fig:example}
\end{figure}

It is evident from these pictures that the pseudo $S$-spectrum of $L_{q}$ behaves differently for different values of $\varepsilon.$ In particular, $\Lambda^{S}_{\varepsilon}(T)$ does not intersect the X-axis (see Figure (C)) for  $\varepsilon<< 0.1$ and on the other hand,  the intersection with $X$-axis is connected (see Figure (B)) for $\varepsilon > 0.1$. This illustrates Theorem \ref{Theorem: emptyorconnected}.
\end{example}

\subsection{$G_1$-class}
 Recall that a bounded linear operator $A$ on a complex Hilbert space  is said to satisfy $G_1$-condition (see  \cite{Putnam, Stampfli}) if the resolvent operator satisfies the first order rate of growth condition with respect to the spectrum:
\begin{equation*}
    \|(zI-A)^{-1}\|=\frac{1}{\text{dist} \big(z,\sigma(A)\big)}\; \text{for all}\; z \in \rho(A), \; \text{the resolvent set of}\; A.
\end{equation*}  
We now propose a quaternionic analogue of $G_{1}$-condition, that is, the first order rate of growth condition for pseudo resolvent operator with respect to the $S$-spectrum. We draw motivation from Lemma \ref{lemma: distance formula} to study the class of operators satisfying  Equation \eqref{Equation: normcondn}. 
\begin{definition}
 An operator $T \in \mathcal{B}(\mathcal{H})$ is said to satisfy $G_{1}$-condition (or to be of $G_{1}$-class) if 
 \begin{align*}
\|Q_{q}(T)^{-1}\|=\frac{1}{\inf\Big\{ |(\mu-q)(\mu-\bar{q})|:\; \mu \in \sigma_{S}(T)\cap \mathbb{C}\Big\}} \, \, \text{for all}\;  q \in \rrho_S(T).
\end{align*}
Equivalently, $T$ is of $G_{1}$-class if and only if $\Lambda^{S}_{\varepsilon}(T) = S^{\mathbb{H}}(\sigma_{S}(T), \varepsilon)$ for every $\varepsilon >0.$
\end{definition}

Clearly, every quaternionic normal operator is of $G_1$-class. However, it may contain non-normal operators. For instance, every subnormal quaternionic operator is of $G_{1}$-class. Let us recall a situation from the theory of complex Hilbert spaces and compare it with the quaternion setup. If $A$ is a complex hyponormal operator, then $(A-zI)$ is hyponormal and hence $(A-zI)^{-1}$ is hyponormal for every $z\in \mathbb{C}\setminus \sigma(A)$ (see \cite{Stampfli}). Whereas, in a quaternionic Hilbert space, the operator $Q_{q}(T)$ may not be hyponormal even when $T$ is hyponormal. As one can see from Problem 209 of   \cite{Halmos} that the square of a hyponormal complex linear operator need not be hyponormal and we provide an example of quaternionic operator by adopting the idea of Halmos to see that $T$ being hyponormal may not guarantee the  hyponormality of  $Q_{q}(T)^{-1}$. 

\begin{example}\label{Example: hypobutnotsquarehypo}
Define $T \colon \ell^{2}(\mathbb{N}, \mathbb{H}) \to \ell^{2}(\mathbb{N}, \mathbb{H})$ by 
\begin{equation*}
    T(x_{1}, x_{2}, x_{3}, \cdots ) = (x_{2}, x_{3}+2x_{1}, x_{4}+2x_{2}, x_{5}+2x_{3}, x_{6}+2x_{4}, \cdots),
 \end{equation*}
 for all $(x_{1}, x_{2}, x_{3}, \cdots) \in \ell^{2}(\mathbb{N}, \mathbb{H}).$ Then $T$ is a bounded quaternionic operator with $\|T\|\leq 3.$ By the usual computation, we get the adjoint $T^{\ast}$ of $T$ as
 \begin{equation*}
     T^{\ast}(x_{1}, x_{2}, x_{3}, \cdots) = (2x_{2}, x_{1}+2x_{3}, x_{2}+2x_{4}, x_{3}+2x_{5}, x_{4}+2x_{6}, x_{5}+2x_{7}, \cdots).
 \end{equation*}
 Then $T$ is hyponormal operator since 
\begin{align*}
    \big\langle (x_{1}, x_{2}, x_{3}, \cdots), \; T^{\ast}T-TT^{\ast}(x_{1}, x_{2}, x_{3}, \cdots) \big\rangle &= \big\langle (x_{1}, x_{2}, x_{3}, \cdots),\; (3x_{1}, 0, 0, \cdots)\big\rangle \\
    &= 3|x_{1}|^{2}\\
     &\geq 0.
 \end{align*} 
Now we shall see that  $T^{2}$ is not hyponormal. The operator $T^{2}$ and its adjoint ${T^{2}}^{*}$ respectively given as follows:
 \begin{equation*}
     T^{2}(x_{1}, x_{2}, x_{3}, \cdots) = (x_{3}+2x_{1},\; x_{4}+4x_{2},\; x_{5}+4x_{3}+4x_{1},\; x_{6}+4x_{4}+4x_{2},\; x_{7}+4x_{5}+4x_{3},\; \cdots)
 \end{equation*}
 and 
 \begin{equation*}
     {T^{2}}^{\ast}(x_{1}, x_{2}, x_{3}, \cdots) = (4x_{3}+2x_{1},\; 4x_{4}+4x_{2},\; 4x_{5}+4x_{3}+x_{1},\; 4x_{6}+4x_{4}+x_{2},\; 4x_{7}+4x_{5}+x_{3},\; \cdots).
 \end{equation*}
 Now by computing ${T^{2}}^{\ast}T^{2}-T^{2}{T^{2}}^{\ast}$ and simplifying, we get 
 \begin{equation*}
     {T^{2}}^{\ast}T^{2}-T^{2}{T^{2}}^{\ast}(x_{1}, x_{2}, x_{3}, \cdots) = (6x_{3}+15x_{1},\; 15x_{2},\; 6x_{1},\; 0,\;0,\; \cdots).
 \end{equation*}
 Therefore, 
 \begin{align*}
    \big\langle (x_{1}, x_{2}, x_{3}, \cdots),\;  &{T^{2}}^{\ast} T^{2}- T^{2}{T^{2}}^{\ast}(x_{1}, x_{2}, x_{3}, \cdots)\big\rangle\\
    &= \big\langle (x_{1}, x_{2}, x_{3}, \cdots),\;  (6x_{3}+15x_{1},\; 15x_{2},\; 6x_{1},\; 0,\;0,\; \cdots)\big\rangle\\
    &= 15 |x_{1}|^{2} + 6 \overline{x}_{1} x_{3}+ 15 |x_{2}|^{2}+ 6\overline{x}_{3}x_{1}. 
 \end{align*}
 If we put $x_{1} = 1, x_{2} = 0, x_{3}= -2$ and $x_{n} = 0 $ for all $n \geq 4$, then value of the above inner product is $-9$. Hence $T^{2}$ is not hyponormal.

\end{example}

 From the learnings of Example \ref{Example: hypobutnotsquarehypo} there is a necessity to introduce a new subclass of hyponormal operator, we call it polynomially hyponormal operators. Though this notion is available in the classical setup, we have discussed here with a suitable modification.   
\begin{definition}
An opertor $T \in \mathcal{B}(\mathcal{H})$ is said to be \textit{ polynomially hyponormal} if $\mathcal{P}(T)$ is hyponormal for every polynomial $\mathcal{P} \in \mathbb{R}[x]$, the class of all real polynomial in the quaternion variable $x$.
\end{definition}

\begin{theorem} \label{Theorem: polynomiallyhypoisG1} Every polynomially hyponormal quaternionic operator is of $G_{1}$-class. 
\end{theorem}
\begin{proof}
Suppose that $T\in \mathcal{B}(\mathcal{H})$ is polynomially hyponormal and  $q\in \rrho_{S}(T)$ then  $Q_{q}(T)$ is hyponormal since $Q_{q}(T) = p(T)$ where $p(x) = x^{2}-2 re(q)x+ |q|^{2}\in \mathbb{R}[x]$. So the inverse $Q_{q}(T)^{-1}$ also hyponormal.  It follows that $\|(Q_{q}(T)^{-1})^{n}\| = \|Q_{q}(T)^{-1}\|^{n}$, for every $n \in \mathbb{N}.$ By the spectral radius formula, we see that  $\|Q_{q}(T)^{-1}\|=r_S(Q_q^{-1})$. Hence $T$ is of $G_{1}$-class.
\end{proof}

The following inclusion holds for quaternionic operators. 
\begin{equation*}
\text{normal }\subseteq \text{ subnormal } \subseteq \text{ polynomially hyponormal} \subseteq G_{1} \text{ -class}.
\end{equation*}

We conclude our discussion with an example of a non-normal quaternionic operator which is of $G_{1}$-class.
\begin{example}
Consider the right shift operator $R: \ell^2(\mathbb{N}, \mathbb{H})\rightarrow \ell^2(\mathbb{N}, \mathbb{H}) $ defined by \begin{equation*}
    R(x_1,x_2,x_3, \cdots)=(0,x_1,x_2,\cdots),\; \text{for all}\;  (x_1,x_2,x_3, \cdots)\in \ell^2(\mathbb{N}, \mathbb{H}).
\end{equation*}
Note that $R$ is not a normal operator. Further, $\sigma_S(R)=\{q\in \mathbb{H}: |q|=1\}.$
In fact, for any $p\in \mathbb{R}[x]$, the operator  $p(R)$ is subnormal. Then $p(R)$ is hyponormal. It follows that $R$ is polynomially hyponormal. Thus by Theorem \ref{Theorem: polynomiallyhypoisG1}, we see that $R$ is of $G_{1}$-class. Therefore,
\begin{equation*}
     \Lambda_\varepsilon^S(R)= S^\mathbb{H}(\sigma_S(R), \varepsilon) = \big\{ q \in \mathbb{H}:\; \inf\limits_{\mu \in \mathbb{C},\;|\mu| =1} |\mu^{2} - 2 \text{Re}(q)\mu+1| \leq \varepsilon\big\}
\end{equation*}
 for all $\varepsilon>0$.
\end{example}

\subsection*{Declaration of competing interest}

The authors declare that there is no competing interest.



\subsection*{Acknowledgment}
We thank the referee for a careful reading of the article and useful suggestions.
Also, we thank the Statistics and Mathematics Unit, Indian Statistical Institute Bangalore Centre for its support and for providing an environment to initiate collaboration. The first named author is supported by the postdoctoral fellowship of Weizmann Institute of Science, Israel. The second named author expresses his gratitude to the  Indian Institute of Science Education and Research (IISER) Mohali, India for providing necessary research facilities to finish this work. 


\end{document}